**Michiel Hazewinkel**
Direct line: +31-20-5924204
Secretary: +31-20-5924089
Fax: +31-20-5924166
E-mail: mich@cwi.nl

CWI
POBox 94079
1090GB Amsterdam

# Word Hopf algebras

by
*Michiel Hazewinkel*
*CWI*
*POBox 94079*
*1090GB Amsterdam*
*The Netherlands*

**Abstract**. Two important generalizations of the Hopf algebra of symmetric functions are the Hopf algebra of noncommutative symmetric functions and its graded dual the Hopf algebra of quasisymmetric functions. A common generalization of the latter is the selfdual Hopf algebra of permutations ($MPR$ Hopf algebra). This latter Hopf algebra can be seen as a Hopf algebra of endomorphisms of a Hopf algebra. That turns out to be a fruitful way of looking at things and gives rise to wide ranging further generalizations such as the word Hopf algebra and the double word Hopf algebra (and many more).

## 1. Introduction and motivation.

It is well known by now that one of the better ways of looking at symmetric functions is to consider an infinity of indeterminates and to look at the basic structures on the symmetric functions as defined by the one of a connected graded selfdual Hopf algebra over the integers. By no means all of the extremely rich structure is captured this way; but a lot is and this makes an extremely good start.

To make this lecture more or less selfcontained a few words about coalgebras and Hopf algebras will be said a bit later in (section 2); mainly in terms of some specific examples.

More or less recently two important generalizations of the symmetric functions have appeared, the Hopf algebra $NSymm$ of noncommutative symmetric functions, and its graded dual the Hopf algebra $QSymm$ of quasisymmetric functions. Many of the structures of the Hopf algebra $Symm$ of symmetric functions have natural analogues for these, see e.g. the two surveys (1, 2), see also the seminal paper (3) and the many subsequent papers by several of the same authors and others. But one thing, besides many others, is lost: they are not selfdual, the one being maximally noncommutative and the other, its dual, commutative. That is one good reason to look for a common generalization that is selfdual.

There is such a generalization, the $MPR$ Hopf algebra of permutations. This is the first of the 'word Hopf algebras' with which this lecture is concerned. It will be described in detail later in section 3. I got interested in it for three reasons. First the formulas defining it seem at first sight rather mysterious and one wonders where they come from; second, one wonders whether there are more Hopf algebras like it; third, for a number of reasons $MPR$ seems a bit small and one wonders if there are more general objects with more room for interesting families of endomorphisms.

In this lecture everything is done over the integers and it is assumed that the underlying Abelian groups are free; everything works also over arbitrary commutative rings $k$ with unit element, in particular fields, with the assumption that the underlying modules of the algebras, coalgebras, bialgebras, Hopf algebras, that occur, are free.

Much more on the subject of this lecture can be found in (4).

## 2. Coalgebras, bialgebras, and Hopf algebras

Here are a few words on coalgebras and Hopf algebras as promised. For much more see e.g. (5, 6, 7, 8).

*Coalgebras.*[1] An algebra over $\mathbf{Z}$ is a $\mathbf{Z}$-module $A$, i.e. an Abelian group, with a bilinear map $A \times A \longrightarrow A$. That is, there is a recipe for taking two elements of $A$ and composing them to get a third new element made up of these two. At this stage we are not concerned with any other properties this composition law might have (such as associativity, Lie, alternative, ...).

Of course a bilinear map $m : A \times A \longrightarrow A$ is the same thing as a $\mathbf{Z}$-module morphism

$$m : A \otimes A \longrightarrow A \tag{2.1}$$

A good example is the word-concatenation-algebra. Let $\mathcal{X}$ be an alphabeth, and let $\mathcal{X}^*$ be the monoid (under concatenation) of all words over $\mathcal{X}$ including the empty word. Let $\mathbf{Z}[\mathcal{X}^*]$ be the free module over $\mathbf{Z}$ with basis $\mathcal{X}^*$. Define a multiplication on $\mathbf{Z}[\mathcal{X}^*]$ by assigning to two words $v, w \in \mathcal{X}^*$ their concatenation (denoted by $*$)

$$m(v,w) = v * w \tag{2.2}$$

and extending bilinearly. Thus if $v = [4,3,5,1]$ and $w = [7,1,1]$ are two words over the natural numbers $\mathbf{N} = \{1,2,\cdots\}$ their concatenation is $v * w = [4,3,5,1,7,1,1]$. This multiplication is associative and the empty word, [ ], serves as a two sided unit element. This algebra is of course also the monoid algebra on $\mathcal{X}^*$ over $\mathbf{Z}$ and the free associative algebra on $\mathcal{X}$ over $\mathbf{Z}$.

Sort of dually, a coalgebra over $\mathbf{Z}$ is a $\mathbf{Z}$-module complete with a decomposition map

$$\mu : C \longrightarrow C \otimes C \tag{2.3}$$

As in the algebra case above, think of the $\mathbf{Z}$-module $C$ as a free module over $\mathbf{Z}$ with as basis some kind of objects (such as words). Then $\mu$ of such a basis object gives all possible ways of decomposing that object into two other basis objects. A good example is the word-cut-coalgebra structure on $\mathbf{Z}[\mathcal{X}^*]$ which is given by

$$\mu([x_1, x_2, \cdots, x_m]) = \sum_{i=0}^{m} [x_1, \cdots, x_i] \otimes [x_{i+1}, \cdots, x_m] \tag{2.4}$$

For example

$$\begin{aligned}\mu([4,3,5,1]) = [\,] \otimes [4,3,5,1] + [4] \otimes [3,5,1] + [4,3] \otimes [5,1] + \\ + [4,3,5] \otimes [1] + [4,3,5,1] \otimes [\,]\end{aligned} \tag{2.5}$$

[1] In a way, historically speaking, coalgebras are older than algebras. They are related to 'addition formuas' such as the formulas of trigonometry like $\sin(x+y) = \sin(x)\cos(y) + \cos(x)\sin(y)$, $\cos(x+y) = \cos(x)\cos(y) - \sin(x)\sin(y)$. See (8) or (13) for how this comes about.

Coassociativity says that the decomposition rule $\mu$ satisfies

$$(id \otimes \mu)\mu = (\mu \otimes id)\mu : C \longrightarrow C \otimes C \otimes C \tag{2.6}$$

which in the present context can be thought of as follows. If one wants to break up an object into three pieces, than first break it up into two pieces (in all possible ways) and than break up the left halfs into two pieces in all possible ways. Or, after the first step, break the right half into two pieces in all possible ways. Coassociativity now says that it does not matter which of the two procedures is followed.

The word-cut-coalgebra just defined is coassociative.

There is also a counit morphism given by

$$\begin{aligned} &\varepsilon([\,]) = 1 \\ &\varepsilon(w) = 0 \ \text{ if } \ w \text{ has length } \ \geq 1 \end{aligned} \tag{2.7}$$

which satisfies $(id \otimes \varepsilon)\mu = id$, $(\varepsilon \otimes id)\mu = id$.

Another good example is $CoF(\mathbf{Z})$ which as an Abelian group is the free group with basis $Z_n$, $n = 0,1,2,\cdots$ . Think of $Z_n$ as representing a horizontal strip consisting of $n$ unit boxes. Such a strip decomposes into two strips of $i$ and $n-i$ boxes respectively, $i = 0,1,\cdots,n$, (including two trivial decompositions).This suggests the decomposition formula

$$\mu(Z_n) = \sum_{i+j=n} Z_i \otimes Z_j \tag{2.8}$$

which is indeed coassociative. There is also a counit morphism $\varepsilon(Z_0) = 1$, $\varepsilon(Z_n) = 0$ for $n > 0$.

*Bialgebras*. A bialgebra $B$ over **Z**, or a **Z**-bialgebra is a **Z**-module $B$ equipped with a multiplication $m$, a comultiplication $\mu$, a unit morphism $e$ (this means that $e(1)$ is the unit element for the multiplication), and a counit morphism $\varepsilon$

$$\begin{aligned} &m: \ B \otimes B \longrightarrow B, \quad e: \ \mathbf{Z} \longrightarrow B \\ &\mu: \ B \longrightarrow B \otimes B, \quad \varepsilon: \ B \longrightarrow \mathbf{Z} \end{aligned} \tag{2.9}$$

such that $(B,m,e)$ is an associative algebra with unit, $(B,\mu,\varepsilon)$ is a coassociative coalgebra with counit, and

$$m \text{ and } e \text{ are coalgebra morphisms} \tag{2.10}$$

$$\mu \text{ and } \varepsilon \text{ are algebra morphisms} \tag{2.11}$$

Here $B \otimes B$ is given the tensor product algebra and coalgebra structures. That is

$$m_{B\otimes B}(a \otimes b \otimes c \otimes d) = ac \otimes bd, \ e_{B\otimes B}(1) = e_B(1) \otimes e_B(1) \tag{2.12}$$

and if

$$\mu(a) = \sum_i a_{i,1} \otimes a_{i,2}, \quad \mu(b) = \sum_j b_{j,1} \otimes b_{j,2}$$

$$\mu_{B \quad B}(a \quad b) = \sum_{i,j} a_{i,1} \quad a_{i,2} \quad b_{j,1} \quad b_{j,2}, \quad \varepsilon_{B \quad B}(a \quad b) = \varepsilon(a)\varepsilon(b) \tag{2.13}$$

And **Z** is given the trivial algebra and coalgebra structures, obtained from the natural isomorphism **Z** **Z** **Z**.

The conditions (1.10) and (1.11) are equivalent and mostly one focusses on (1.11). This compatibility condition, often called the Hopf condition, is far from saying that comultiplication (=decomposition) undoes multiplication (=composition). Rather it says that the two are sort of orthogonal to each other.

Thus $\mathbf{Z}[\mathcal{X}\ ]$ with concatenation, see (1.2), as multiplication, and cut, see (1.4), as comultiplication, is *not* a bialgebra.

Think of the comultiplication as giving all ways of composing an object into a left part and a right part and of multiplication as vertical composition. Then in words the Hopf condition says that all left parts of the decomposition of a product are all products of left parts of the two factors and the corresponding right parts are the products of the corresponding right parts of the decompositions of the factors involved. In pictures (stacking things vertically is composition; cutting things into left and right parts is decomposition):

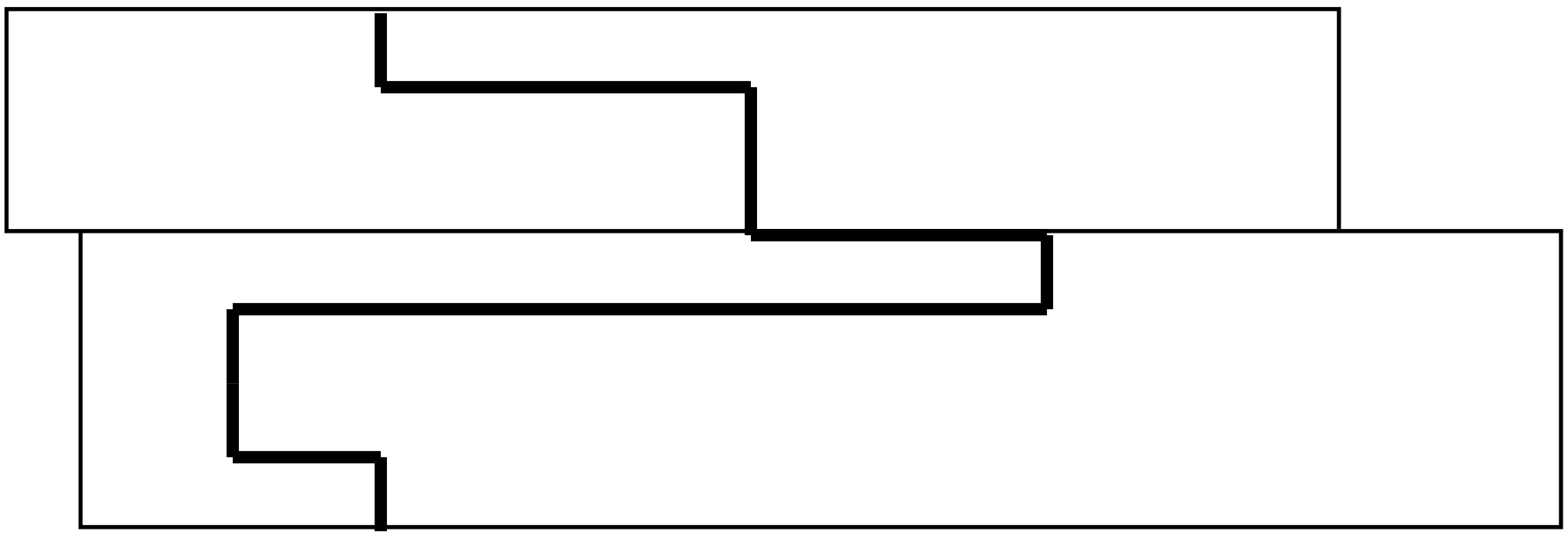

Here the top box is object 1 and the bottom one is object 2. The two boxes together are the composition of object 1 and object 2 (in that order). The first three segments of the jagged heavy line give a decomposition of object 1 into a left part L1 and a right part R1 and the last five segments of that heavy jagged line give a decomposition into a left part L2 and a right part R2 of object 2. The complete heavy broken line gives the corresponding decomposition of the product of object 1 and 2 into a left part that is the product of L1 and L2 and a right part that is the product of R1 and R2.

In fact there is a most important and beautiful bialgebra (which is in fact a Hopf algebra) that looks almost exactly like this.

Consider stacks of rows of unit boxes

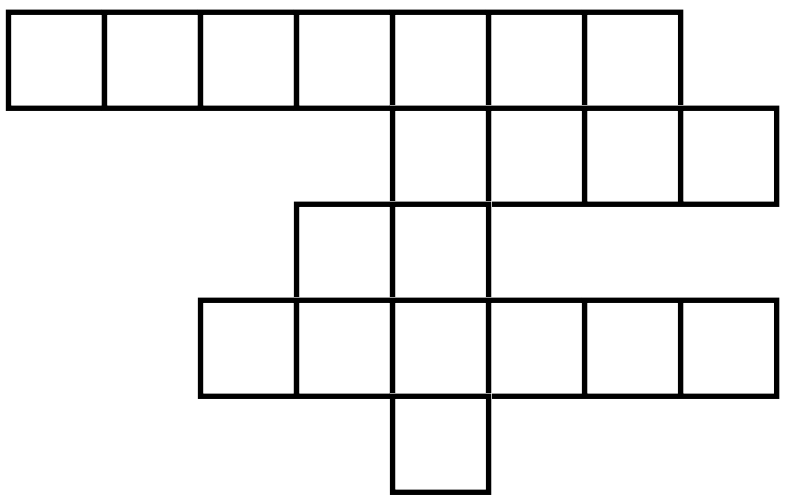

Here two such stacks are considered equivalent if they have the same number of boxes in each layer. Thus it only matters how many boxes there are in each layer. The case depicted is hence given by the word [7,4,2,6,1] over the positive integers $\mathbf{N}$. The empty word is permitted and corresponds to any stack of empty layers. There may be empty layers; these are ignored. The possible decompostions of a stack are obtained by cutting each layer into two parts. Two of these for the example at hand are indicated below.

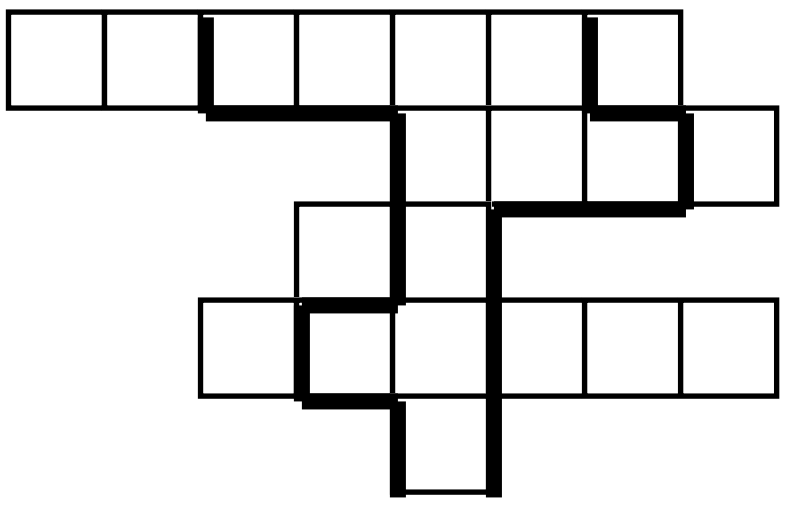

These two correspond to the decompositions [2,1,1]   [5,4,1,5,1] and [6,3,2,3,1]   [1,1,3]. A convenient way of encoding this algebraically is to consider the free associative algebra $\mathbf{Z}\ Z$ over the integers in the indeterminates $Z_1, Z_2, Z_3, \cdots$. A basis (as a free Abelian group) for this algebra is given by 1 and all (noncommutative) monomials $Z_{i_1} Z_{i_2} \cdots Z_{i_m}$. This monomial encodes the stack with layers of $i_1, i_2, \cdots, i_m$ boxes. Thus the example above corresponds to the monomial $Z_7 Z_4 Z_2 Z_6 Z_1$.

The comultiplication on $\mathbf{Z}\ Z$ is given by the algebra homomorphism determined by

$$\mu(Z_n) = \sum_{i+j=n} Z_i \quad Z_j \quad \text{where } Z_0 = 1 \text{ and } i,j \quad \mathbf{N} \quad \{0\} = \{0,1,2,3,\cdots\} \tag{2.14}$$

This, together with the counit

$$\varepsilon(Z_n) = 0, \ n \geq 1 \tag{2.15}$$

defines the bialgebra (Hopf algebra) *NSymm* of noncommutative symmetric functions. Note that there is no difficulty about verifying the Hopf property because $\mathbf{Z}\ Z$ is a free algebra and the comultiplication is defined as the algebra morphism determined by (2.14).

*Hopf algebras*. A Hopf algebra $H$ is a bialgebra with one additional piece of structure, a morphism of modules $\iota : H \longrightarrow H$ that is required to satisfy

$$\text{conv}(\text{id},\iota) = \text{conv}(\iota,\text{id}) = \varepsilon \circ e : \ H \longrightarrow H \tag{2.16}$$

Here $\text{conv}(f,g)$ is the convolution product of two module endomorphisms $f,g: \ H \longrightarrow H$ which is defined as the composite

$$H \longrightarrow H \quad H \longrightarrow H \quad H \longrightarrow H \tag{2.17}$$

The module endomorphism $\varepsilon \circ e: \ H \longrightarrow H$ is the identity for the convolution product.

If one exists, an antipode is unique. In many cases the existence of an antipode follows from other properties of a bialgebra. Antipodes paly a role in Hopf algebra theory that is somewhat analogous to that of inverses in group theory.

*Example*. Monoid bialgebra of a finite monoid. Let $G$ be a finite monoid. Let $kG$ be the monoid algebra. Define

$\mu: kG \longrightarrow kG \otimes kG$ and $\varepsilon: kG \longrightarrow k$ by

$$\mu(g) = g \otimes g, \quad \varepsilon(g) = 1, \quad g \in G \tag{2.18}$$

It is easy to check that these are algebra morphisms. Indeed it suffices to check this on the basis $\{g: g \in G\}$ and then $\mu(gh) = gh \otimes gh = (g \otimes g)(h \otimes h) = \mu(g)\mu(h)$, $\varepsilon(gh) = 1 = \varepsilon(g)\varepsilon(h)$. If $G$ is a group, there is an antipode, viz. $g \mapsto g^{-1}$. In fact a monoid bialgebra has an antipode if and only if the defining monoid is a group.

*Example*. The shuffle Hopf algebra.
Consider the monoid $\mathbf{N}^*$ all words over the alphabeth $\mathbf{N}$ of positive integers (including the empty word) and let *Shuffle* be the free $\mathbf{Z}$-module (Abelian group) with basis $\mathbf{N}^*$. The shuffle product of two words $\alpha = [a_1, a_2, \cdots, a_m]$ and $\beta = [b_1, b_2, \cdots, b_n]$ is defined as follows. Take a 'sofar empty' word with $n+m$ slots. Choose $m$ of the available $n+m$ slots and place in it the natural numbers from $\alpha$ in their original order; place the entries from $\beta$ in their original order in the remaining $n$ slots. The product of the two words $\alpha$ and $\beta$ is the sum (with multiplicities) of all words that can be so obtained. So, for instance

$$[a,b] \times_{sh} [c,d] = [a,b,c,d] + [a,c,b,d] + [a,c,d,b] + [c,a,b,d] + [c,a,d,b] + [c,d,a,b]$$

$$[1] \times_{sh} [1,1,1] = 4[1,1,1,1]$$

This defines a commutative associative multiplication on *Shuffle* for which the empty word is a unit element. Moreover with cut as a comultiplication

$$\mu([a_1, a_2. \cdots, a_m]) = \sum_{i=0}^{m} [a_1, \cdots, a_i] \otimes [a_{i+1}, \cdots, a_m] \tag{2.19}$$

counit $\varepsilon([]) = 1$, $\varepsilon(\alpha) = 0$ if $\lg(\alpha) \geq 1$, where the length of a word $\alpha = [a_1, \cdots, a_m]$ is $\lg(\alpha) = m$, and antipode

$$\iota([a_1, a_2, \cdots, a_m] = (-1)^m [a_m, a_{m-1}, \cdots, a_1]$$

*Shuffle* becomes a Hopf algebra. It is a not entirely trivial (and good) exercise to check the Hopf property for this Hopf algebra.

*Grading.* A graded module over $\mathbf{Z}$ is a module $M$ that comes with a direct sum decomposition into submodules

$$M = \bigoplus_{i=0}^{\infty} M_i \tag{2.20}$$

The elements of the direct summand $M_i$ are called homogeneous of degree $i$. A morphism $\varphi$ between graded $\mathbf{Z}$-modules $M$ and $N$ is called homogenous if $\varphi(M_i) \subseteq N_i$ for all $i$.

A graded module $M$ is connected if $M_0 = \mathbf{Z}$.[2]

The tensor product (over $\mathbf{Z}$) of two graded modules $M$ and $N$ is made a graded module by

[2] This terminology comes from algeebraic topology. The cohomology and homology of connected spaces are connected graded modules in this sense.

$$(M \otimes N)_n = \bigoplus_{i+j=n} M_i \otimes N_j \tag{2.21}$$

An algebra $A$ whose underlying module is graded is a graded algebra if the multiplication morphism $m\colon A \otimes A \longrightarrow A$ is homogeneous (and $e(\mathbf{Z}) \subset A_0$). A coalgebra $C$ whose underlying module is graded is a graded coalgebra if the comultiplication morphism $\mu\colon C \longrightarrow C \otimes C$ is homogeous (and $\varepsilon(C_i) = 0$ for all $i > 0$).

A bialgebra $H$ is a graded bialgebra if it is a both a graded algebra and a graded coalgebra. Both bialgebras $NSymm$ and $Shuffle$ are connected graded, with the grading given, respectively, by

$$\deg(Z_{a_1} Z_{a_2} \cdots Z_{a_m}) = a_1 + a_2 + \cdots + a_m \tag{2.22}$$

$$\deg([a_1, a_2, \cdots, a_m]) = a_1 + a_2 + \cdots + a_m \tag{2.23}$$

For a connected graded bialgebra the existence of an antipode is automatic and is homogeneous making it into a graded Hopf algebra.

## 3. MPR, the Hopf algebra of permutations.

The original motivation for the studies on which this lecture (partially) reports[3] comes from the $MPR$ Hopf algebra (Hopf algebra of permutations) introduced and studied by Malvenuto, Poirier and Reutenauer, (9, 10). Here is a description.

As a graded Abelian group $MPR$ has a basis consisting of the empty permutation and all permutations on $n$ letters, $n = 1, 2, \cdots$ . Thus

$$MPR = \mathbf{Z} \oplus \bigoplus_{n=1}^{\infty} \mathbf{Z} S_n \tag{3.1}$$

where $S_n$ is the symmetric group on $n$ letters. $\mathbf{Z} S_n$ is the homogeneous summand of degree $n$.

Permutations on $m$ letters are written as words, with a word $\alpha = [a_1, a_2, \cdots, a_m]$ on the alphabet $\{1, \cdots, m\}$ which has no repeats, corresponding to the permutation $i \mapsto a_i$, $i = 1, \cdots, m$. The empty permutation corresponds to the empty word [ ], which is the basis element of the first summand in (3.1).

*Multiplication on MPR*. With these notations the multiplication on $MPR$ can be described as follows. The empty word serves as the unit element, and if $\alpha = [a_1, \cdots, a_m]$ and $\beta = [b_1, \cdots, b_n]$ are two permutation words their product is

$$m(\alpha \otimes \beta) = [a_1, \cdots a_m] \times_{sh} [m + b_1, \cdots, m + b_n] \tag{3.2}$$

where $\times_{sh}$ stands for the shuffle product of two words, described above. Thus for instance $m_{MPR}([1] \otimes [3,2,1]) = [1] \times_{sh} [4,3,2] = [1,4,3,2] + [4,1,3,2] + [4,3,1,2] + [4,3,2,1]$.

*Standardization*. To describe the comultiplication the notion of standardization is needed. For any word $\alpha = [a_1, a_2, \cdots, a_m]$ without repeats over the alphabet of natural numbers $\mathbf{N} = \{1, 2, 3, \cdots\}$ its standardization is the permutation word $\mathrm{st}(\alpha) = [\varphi(a_1), \cdots, \varphi(a_m)]$ where $\varphi : \{a_1, \cdots, a_m\} \longrightarrow \{1, \cdots, m\}$ is the unique strictly monotone map between these ordered sets. For instance $\mathrm{st}([5,2,1,8] = [3,2,1,4]$. This notion of standardization is a special case of standardization of words as introduced by Schensted, (11), which applies to all words over $\mathbf{N}$, not only to words without repeats.

[3] Much more can be found in (4).

*Comultiplication on MPR*. The comultiplication on $MPR$ is now defined by

$$\mu(\alpha) = \sum_{\alpha' \quad \alpha'' = \alpha} st(\alpha') \quad st(\alpha'') \tag{3.3}$$

where the sum is over all cuts of $\alpha$ , that is all pairs of words $(\alpha', \alpha'')$ whose concatenation $\alpha' \quad \alpha''$ is equal to $\alpha$.

With $1 = [\ ]$ as unit and a counit $\varepsilon$ defined by $\varepsilon(\mathbf{Z}S_n) = 0$ and $\varepsilon([\ ]) = 1$, $MPR$ becomes a bialgebra; there is also (of course, given the graded setting) an antipode, making $MPR$ a Hopf algebra.

The question which intrigued me is now: How does one dream up such formulas and are there natural generalizations?

For instance, there is a very natural generalization of (3.2) to arbitrary words as follows. Define the height of a word $\alpha = [a_1, a_2, \cdots, a_m]$ as $\mathrm{ht}(\alpha) = \max\{a_1, \cdots, a_m\}$ and define the product of two arbitrary words $\alpha = [a_1, a_2, \cdots, a_m]$ and $\beta = [b_1, \cdots, b_n]$ as

$$m_{WHA}(\alpha \quad \beta) = [a_1, \cdots, a_m] \times_{sh} [\mathrm{ht}(\alpha) + b_1, \cdots, \mathrm{ht}(\alpha) + b_n] \tag{3.4}$$

(Note that this agrees with (3.2) when $\alpha$ is a permutation word.) The question now arises whether there is a corresponding comultiplication turning the graded Abelian group with as basis all words over the natural numbers into a Hopf algebra (with (3.4) as its multiplication). As it turns out there is and the result is what I call the word Hopf algebra (*WHA*). Here are some examples of the comultiplication

$$\alpha = [3,2,7,2,4],\ \mu(\alpha) = 1 \quad \alpha + [1] \quad [2,6,2,3] + [3,2,6,2] \quad [1] + \alpha \quad 1$$

$$\alpha = [7,3,2,2,4],$$
$$\mu(\alpha) = 1 \quad \alpha + [3] \quad [3,2,2,4] + [4,1] \quad [2,2,3] + [6,3,2,2] \quad [1] + \alpha \quad 1$$

and it does not seem easy (to me) to guess at these formulas, or more precisely at the general recipe behind them.

There are also many more generalizations such as *dWHA*, the double word Hopf algebra.

There is even more structure to $MPR$. For one thing there is a grading given by $\deg([s_1, s_2, \cdots, s_m]) = m$ making it a graded connected Hopf algebra. There is also a second multiplication given by composition of permutations if these are of equal length and zero otherwise.
Further there is a nondegenerate (but non positive definite) inner product

$$\langle \sigma, \tau \rangle = \begin{cases} 1 & \text{if } \tau = \sigma^{-1} \\ 0 & \text{otherwise} \end{cases} \tag{3.5}$$

and $MPR$ is selfdual with respect to this inner product, meaning that

$$\langle m(\rho \quad \sigma), \tau \rangle = \langle \rho \quad \sigma, \mu(\tau) \rangle \tag{3.6}$$

where on the right hand side of (3.6) the inner product is the tensor inner product given by

$$\langle \rho \otimes \sigma, \rho' \otimes \sigma' \rangle = \langle \rho, \rho' \rangle \langle \sigma, \sigma' \rangle .$$

And so supplementary questions arise: where do all these come from and do they generalize?

## 4. Seemingly irrelevant intermezzo: End($H$).

In this section $H$ is a finite rank Hopf algebra (over $\mathbf{Z}$), meaning that the underlying Abelian group is free of finite rank. Let $\mathrm{End}(H)$ denote the Abelian group of Abelian group endomorphisms of $H$.

For any Abelian group $H$ there are natural morphisms

$$\begin{aligned} &H \otimes H^* \longrightarrow \mathrm{End}(H), \quad u \otimes f \mapsto \varphi, \ \varphi(v) = f(v)u \\ &\mathrm{End}(H) \otimes \mathrm{End}(H) \longrightarrow \mathrm{End}(H \otimes H), \ (f \otimes g)(u \otimes v) = f(u) \otimes g(v) \end{aligned} \tag{4.1}$$

If $H$ is a free Abelian group of finite rank these two Abelian group endomorphisms are isomorphisms. (But if $H$ is of infinite rank this is definitely not the case.) Here $H^* = Hom(H, \mathbf{Z})$ the linear dual of $H$, which is a Hopf algebra with comultiplication $m^*$, multiplication $\mu^*$, unit $\varepsilon^*$, counit $e^*$ and antipode $\iota^*$.

Now, as the tensor product of two Hopf algebras $H \otimes H^*$ carries a Hopf algebra structure. Here is how that structure looks in terms of endomorphisms, i.e on $\mathrm{End}(H)$.

*Convolution*. The multiplication is convolution. Let $f$ and $g$ be two elements of $\mathrm{End}(H)$. Then their product is the composite

$$H \xrightarrow{\ \mu\ } H \otimes H \xrightarrow{\ f \otimes g\ } H \otimes H \xrightarrow{\ m\ } H \tag{4.2}$$

*Coconvolution*. The comultiplication is given as follows. Let $f$ be an element of $\mathrm{End}(H)$. Then $\mu(f)$ is the following element of $\mathrm{End}(H) \otimes \mathrm{End}(H)$

$$H \otimes H \xrightarrow{\ m\ } H \xrightarrow{\ f\ } H \xrightarrow{\ \mu\ } H \otimes H \tag{4.3}$$

Of course (4.3) defines an element of $\mathrm{End}(H \otimes H)$; but then using the second of the isomorphisms of (4.1) this yields an element of $\mathrm{End}(H) \otimes \mathrm{End}(H)$.

*Unit*. The unit of $\mathrm{End}(H)$ is the endomorphism $H \xrightarrow{\ \varepsilon\ } \mathbf{Z} \xrightarrow{\ e\ } H$ of $H$.

*Counit*. the counit of $\mathrm{End}(H)$ takes an endomorphism $f$ in $\mathrm{End}(H)$ to the number $(\varepsilon \circ f \circ e)(1)$.

*Antipode*. Finally the antipode of $\mathrm{End}(H)$ is given by $f \mapsto \iota \circ f \circ \iota$.

The Hopf algebra $\mathrm{End}(H)$ is also selfdual. This is easiest seen in its guise $H \otimes H^*$. Indeed

$$(H \otimes H^*)^* \longrightarrow H^* \otimes H^{**} \longleftarrow H^* \otimes H \longrightarrow H \otimes H^*$$

The corresponding inner product is

$$(H \otimes H^*) \times (H \otimes H^*) \longrightarrow \mathbf{Z}, \ \langle u \otimes f, v \otimes g \rangle = f(v)g(u)$$

which is nondegenerate but not positive definite.

At the level of $\mathrm{End}(H)$ the autoduality is a combination of a canonical pairing

$$\mathrm{End}(H) \quad \mathrm{End}(H\ ) \longrightarrow \mathbf{Z} \tag{4.4}$$

and the isomorphism $\mathrm{End}(H) \longrightarrow \mathrm{End}(H\ )$ that assigns to an endomorphism of $H$ the dual endomorphism which is an endomorphism of $H$ . The canonical pairing (4.4) is defined as follows. Take a basis $u_1, \cdots, u_n$ of $H$ and let $v^1, \cdots, v^n$ be the dual basis of $H$ (so that $v^i(u_j) = \delta^i_j$ (Kronecker delta). Then

$$\gamma = \sum_{i=1}^{n} u_i \quad v^i \qquad H \quad H \tag{4.5}$$

is a special element of $H \quad H$ that is independent of the choice of basis. The easiest way to see this is to remark that the element $\gamma$ is the image of $1 \quad \mathbf{Z}$ under the dual of the evaluation morphism

$$\mathrm{ev}:\ H \quad H \longrightarrow \mathbf{Z}, \quad u \quad f \mapsto f(u) \tag{4.6}$$

The pairing (4.4) is now defined by

$$f \quad g \mapsto (f \quad g)(\gamma) = \sum_{i=1}^{n} f(u_i) g(v^i) \tag{4.7}$$

The elements of $\mathrm{End}(H)$ are endomorphisms. So there is a second multiplication on $\mathrm{End}(H)$, viz composition. This second multiplication is not necessarily distributive over the first (which would make $\mathrm{End}(H)$ a ring object in the category of coalgebras). Still an extra multiplication, i.e. a second way of producing a new element from two given ones, can be very useful even when it has no particular compatibility properties with respect to the other structure present. Examples of this are the many extra ‘multiplications’ of divided power sequences as they are used in (12) to describe a basis over the integers of the Lie algebra of primitives $\mathrm{Prim}(NSymm)$.

Dually there is also a second comultiplication (cocomposition).

4.8. *Open problem*. Which Hopf algebras, necessarily of square rank, are of the form $\mathrm{End}(H)$?

More generally there is a Hopf algebra structure on $\mathrm{Mod}(H, K)$ where $H$ and $K$ are two possibly different Hopf algebras with free finite rank underlying Abelian groups and there is the same open problem vis à vis these Hopf algebras. The same open problems can be considered over other base rings, for instance fields.

Note that practically everything above works perfectly fine for graded Hopf algebras (with graded dual replacing dual). The sole exception is the comultiplication (given by coconvolution). But that is is formidable exception and obstruction.

The seeming irrelevancy of these considerations lies in the following. Of course for a graded Hopf algebra with free finite rank homogeneous summands $H \quad H^{\mathrm{gr}}$ is again a Hopf algebra; but it is a very small part of the Abelian group of homogenous endomorphisms of $H$. Indeed it is easy to check that a homogenous endomorphism of $H$ is in (the image of) $H \quad H^{\mathrm{gr}}$ (in $\mathrm{End}(H)$) iff it has a finite rank image.[4]

Now, as remarked before, the elements of $MPR$ (i.e. permutations) are to be interpreted as homogeneous endomorphisms of the Hopf algebra *Shuffle*. However with the exception of the scalar multiples of the empty permutation none of these endomorphisms lies in $Shuffle \otimes Shuffle^{gr}$.

## 5. $MPR$ as a Hopf algebra of endomorphisms.

The first thing to do is to interprete elements of $MPR$ as endomorphisms (of Abelian groups) of *Shuffle*.

*Permutations as Shuffle endomorphisms*. Let $\sigma = [s_1, \cdots, s_m]$ be a permutation word. The corresponding permutation is of course $j \mapsto s_j$. Further let $\alpha$ be a basis element of *Shuffle*, i.e. a word. Then $\sigma$ takes $\alpha$ to zero unless $\lg(\alpha) = m = \lg(\sigma)$ and if $\alpha = [a_1, \cdots, a_m]$ is of length $m$ then

$$\sigma(\alpha) = [a_{\sigma(1)}, \cdots, a_{\sigma(m)}] = [a_{s_1}, \cdots, a_{s_m}] \qquad (5.1)$$

Note that unless $\sigma$ is the empty permutation (which takes [ ] to [ ] and is zero on all other words) $\sigma$ is never in $Shuffle \otimes Shuffle^{gr}$ (because for every length $\geq 1$ there are infinitely many words of that length).

*Convolution*. Here is how convolution of two permutations works out with this interpretation of permutations as endomorphisms of *Shuffle*. So let $\sigma = [s_1, \cdots, s_m]$ and $\tau = [t_1, \cdots, t_n]$ be two permutation words. Their convolution[5] is given by

$$Shuffle \xrightarrow{\mu_{Sh}} Shuffle \otimes Shuffle \xrightarrow{\sigma \otimes \tau} Shuffle \otimes Shuffle \xrightarrow{m_{Sh}} Shuffle$$

Because $\mu_{Sh}$ is cut, the middle morphism is zero on all terms of $\mu_{Sh}(\alpha)$ for all words $\alpha$ which are not of length $m+n$. And if $\alpha = [a_1, \cdots, a_m, a_{m+1}, \cdots, a_{m+n}]$, $\sigma \otimes \tau$ is zero on all terms of $\mu_{Sh}(\alpha)$ except the summand $[a_1, \cdots, a_m] \otimes [a_{m+1}, \cdots, a_{m+n}]$ and this summand is taken to $[a_{s_1}, \cdots, a_{s_m}] \otimes [a_{m+t_1}, \cdots, a_{m+t_n}]$. And thus the convolution product of $\sigma$ and $\tau$ takes $\alpha$ to the shuffle product $[a_{s_1}, \cdots, a_{s_m}] \times_{Sh} [a_{m+t_1}, \cdots, a_{m+t_n}]$ and it follows immediately that the convolution of $\sigma$ and $\tau$ is equal to to the sum of permutation words $[s_1, \cdots, s_m] \times_{sh} [m+t_1, \cdots m+t_n]$ which is the multiplication on $MPR$ (see section 3 above). There is nothing new about this; this is the way the multiplication of $MPR$ was introduced in (10).

*Coconvolution*. The coconvolution construction takes an endomorphism $f$ of *Shuffle* to the composite morphism

$$Shuffle \otimes Shuffle \xrightarrow{m_{Sh}} Shuffle \xrightarrow{f} Shuffle \xrightarrow{\mu_{Sh}} Shuffle \otimes Shuffle$$

that is, it defines a morphism of Abelian groups

[4] More generally for an infinite rank free Abelian group $M$ or vectorspace an endomorphism is in $M \otimes M$ iff the images of both the morphism itself and its dual are of finite rank.

[5] Quite generally if $C$ is a coalgebra and $A$ is an algebra, the convolution of two morphisms $f, g : C \longrightarrow A$ is the morphism $C \xrightarrow{\mu_C} C \otimes C \xrightarrow{f \otimes g} A \otimes A \xrightarrow{m_A} A$. For the right Hopf algebras this is indeed the classical convolution of functions, see (8) or (7).

$$MPR \subset \mathrm{End}(Shuffle) \xrightarrow{coconv} \mathrm{End}(Shuffle \otimes Shuffle) \tag{5.2}$$

And to turn this into a comultiplication on $MPR$ some sort of projection is needed from the image of coconv in $\mathrm{End}(Shuffle \otimes Shuffle)$ to $MPR \otimes MPR \subset \mathrm{End}(Shuffle) \otimes \mathrm{End}(Shuffle)$. Here is an example how that might be accomplished.

Take $\sigma = [3,1,4,5,2]$ and lets see what the coconvolution of this does to an element of the form $[a_1,a_2] \otimes [b_3,b_4,b_5]$ from $Shuffle \otimes Shuffle$. (Different letters are used for the two copies of $Shuffle$ to make it easier to identify where each is from.)

$$\begin{aligned}
&[a_1,a_2] \otimes [b_3,b_4,b_5] \overset{\times_{sh}}{\longmapsto} [a_1,a_2,b_3,b_4,b_5] + [a_1,b_3,a_2,b_4,b_5] + \\
&+[a_1,b_3,b_4,a_2,b_5] + [a_1,b_3,b_4,b_5,a_2] + [b_3,a_1,a_2,b_4,b_5] + [b_3,a_1,b_4,a_2,b_5] + \\
&+[b_3,a_1,b_4,b_5,a_2] + [b_3,b_4,a_1,a_2,b_5] + [b_3,b_4,a_1,b_5,a_2] + [b_3,b_4,b_5,a_1,a_2] \\
&\underset{\sigma}{\longmapsto} [b_3,a_1,b_4,b_5,a_2] + [a_2,a_1,b_4,b_5,b_3] + \cdots + [b_5,b_3,a_1,a_2,b_4]
\end{aligned}$$

and now cut has to be applied to the final ten terms. But we are looking for an endomorphism of the form $MPR \otimes MPR$. So the only cuts that can contribute are ones which have only $a$'s on the left and only $b$'s on the right. The only one of the final ten terms for which this is possible is the second one and there is (of course) only one cut of this term which qualifies, yielding

$$[a_2,a_1] \otimes [b_4,b_5,b_3] = (\mathrm{st}([3,1]) \otimes \mathrm{st}([4,5,2]))([a_1,a_2] \otimes [b_3,b_4,b_5])$$

And thus according to this procedure the $(\lg = 2) \otimes (\lg = 3)$ component of $\mu_{MPR}([3,1,4,5,2])$ is equal to $\mathrm{st}([3,1]) \otimes \mathrm{st}([4,5,2])$ exactly as it should be according to the description given in section 3 above.

This works in general and gives the correct description of the comultiplication on $MPR$.

All the same this is a very dodgy procedure. For instance, there is a very nice generalization of the way $MPR$ acts as endomorphisms on $Shuffle$ for arbitrary words. This goes as follows. Let now $\sigma = [s_1,\cdots s_n]$ be an arbitrary word of height $m = \mathrm{ht}(\sigma) = \max\{s_1,\cdots,s_n\}$. Then $\sigma$ acts as zero on all basis words $\alpha$ of $Shuffle$ of length $\neq m$ and if $\alpha = [a_1,\cdots a_m]$ is a word of length $m$ $\sigma(\alpha) = [a_{s_1},\cdots,a_{s_n}]$. For permutation words this is the same as (5.1).

Now apply the procedure used above for $MPR$. This does define a comultiplication, but the Hopf property (viz. multiplication is a coalgebra morphism, or, equivalently, the comultiplication is an algebra morphism [6]) fails completely.[7]

However, as it turns out there is a rather different way to make arbitrary words act on $Shuffle$ which does yield a Hopf algebra and for which the dodgy procedure outlined above in the case of $MPR$ does work.

What is needed to use coconvolution to define a Hopf algebra on a suitable module of endomorphisms $E(H)$ is some kind of suitable projection $\mathrm{End}(H \otimes H) \rightarrow E(H) \otimes E(H)$ (where $\mathrm{End}(H) \otimes \mathrm{End}(H)$ is seen as a submodule of $\mathrm{End}(H \otimes H)$ in the natural way). It is largely an open problem what conditions are needed to make this work; see (4) for a preliminary analysis.

[6] This is most always the hardest property to check.

[7] Some two years ago I wasted a couple of months research time to try to fix things up; typical Ptolemaic-epicycle-type thinking, and thoroughly useless.

## 6. *dWHA*, the double word Hopf algebra

Let $\mathcal{X} = \{x_1, x_2, \cdots\}$ be an auxiliary alphabet. A basis as an Abelian group of *dWHA* is formed by pairs of words in the auxiliary alphabet $\mathcal{X}$ with equal support

$$p = \begin{pmatrix} \rho \\ \sigma \end{pmatrix}, \quad \mathrm{supp}(\rho) = \mathrm{supp}(\sigma)$$

The support of a word is the set of different letters that occur in it. I call such pairs of words substitutions.

Here the actual symbols that occur are not important; it is only the patterns of $\rho$ and $\sigma$ relative to each other that are relevant. Thus for instance

$$\begin{pmatrix} [x_1, x_2, x_1, x_3, x_3, x_1, x_4] \\ [x_2, x_3, x_2, x_4, x_1] \end{pmatrix}, \begin{pmatrix} [x_7, x_6, x_7, x_2, x_2, x_7, x_5] \\ [x_6, x_2, x_6, x_5, x_7] \end{pmatrix}, \begin{pmatrix} [y_3, z_4, y_3, x_2, x_2, y_3, x_1] \\ [z_4, x_2, z_4, x_1, y_3] \end{pmatrix} \tag{6.1}$$

all denote the same basis element of *dWHA*. Abstractly

$$p = \begin{pmatrix} \rho \\ \sigma \end{pmatrix}, \ \mathrm{supp}(\rho) = \mathrm{supp}(\sigma) \ \text{ and } \ p = \begin{pmatrix} \rho' \\ \sigma' \end{pmatrix}, \ \mathrm{supp}(\rho') = \mathrm{supp}(\sigma')$$

are the same substitution if and only if there is a bijection $\vartheta : \mathrm{supp}(\rho) \longrightarrow \mathrm{supp}(\rho')$ such that $\vartheta_*(\rho) = \rho'$ and $\vartheta_*(\sigma) = \sigma'$ where $\vartheta_*$ is the map on words induced by $\vartheta$ obtained by applying it to the letters making up a word.

These elements $p$ can (and often should) be thought of as defining endomorphisms of *Shuffle*; more precisely they are recipes for endomorphisms as follows. The 'substitution' $p$ acts as zero on all words over $\mathbf{N}$, i.e. the canonical basis elements of *Shuffle*, that are not of the same pattern as $\rho$ and if a word is of the same pattern as $\rho$ then it is taken into the corresponding basis element of *Shuffle* represented by the pattern $\sigma$.

Thus, for example, if $p$ is the substitution (6.1) and $\alpha = [a_1, a_2, \cdots, a_m]$

$$p(\alpha) = \begin{cases} 0 \ \text{ unless } \ \lg(\alpha) = m = 7 \ \text{ and } \ a_1 = a_3 = a_6, \ a_4 = a_5 \\ [a_2, a_3, a_2, a_4, a_1] \ \text{ if } \ \lg(\alpha) = m = 7 \ \text{ and } \ a_1 = a_3 = a_6, \ a_4 = a_5 \end{cases}$$

Obviously these endomorphisms satisfy a 'homogeneity' property; they act the same everywhere. For instance if $\phi : \mathbf{N} \longrightarrow \mathbf{N}$ is any map and $\phi_* : \mathbf{N} \longrightarrow \mathbf{N}$ denotes the corresponding induced map on words $\phi_*(\alpha) = [\phi(a_1), \phi(a_2), \cdots, \phi(a_m)]$

$$\phi_* \circ p = p \circ \phi_* \tag{6.2}$$

*Open problem.* Characterize the 'recipe endomorphisms' $p$ more precisely.

The next step is to describe the graded Hopf algebra structure on *dWHA*.

*Underlying Abelian group.* The underlying Abelian group is the countable free Abelian group with as basis all substitutions $p$. Included is the 'empty substitution'

$$\begin{pmatrix} [\ ] \\ [\ ] \end{pmatrix}$$

which acts on *Shuffle* by taking the empty word [ ] to itsef and every other basis element of *Shuffle* to zero.

*Grading.* The grading on $dWHA$ is given by

$$\deg(p) = \#\mathrm{supp}(\rho) \tag{6.3}$$

For example the degree of the basis element (6.1) is 4. The degree of the empty substitution is zero and that is the only basis element of degree zero so that the graded Abelian group $dWHA$ is connected. Note that the rank of each homogeneous piece of $dWHA$ is infinite.

*Multiplication.* Let

$$p = \begin{pmatrix} \rho \\ \sigma \end{pmatrix}, \quad p' = \begin{pmatrix} \rho' \\ \sigma' \end{pmatrix} \tag{6.4}$$

be two substitutions. If necessary, first rewrite the second one (or the first one, or both) so that $\mathrm{supp}(\rho) \cap \mathrm{supp}(\rho') = \emptyset$. Then the product of the two substitutions (6.4) is the sum of substitutions

$$m_{dWHA}(p \otimes p') = \begin{pmatrix} \rho \cdot \rho' \\ \sigma \times_{sh} \sigma' \end{pmatrix} \tag{6.5}$$

where $\cdot$ denotes concatenation, $\times_{sh}$ is the shuffle product and if $u = \sigma_1 + \sigma_2 + \cdots + \sigma_r$ is a sum of words with $\mathrm{supp}(\sigma_1) = \mathrm{supp}(\sigma_2) = \cdots = \mathrm{supp}(\sigma_r) = \mathrm{supp}(\rho)$

$$\begin{pmatrix} \rho \\ u \end{pmatrix} = \begin{pmatrix} \rho \\ \sigma_1 \end{pmatrix} + \begin{pmatrix} \rho \\ \sigma_2 \end{pmatrix} + \cdots + \begin{pmatrix} \rho \\ \sigma_r \end{pmatrix}$$

*Unit element.* The unit element is the empty substitution.

It is easy to see that the multiplication is associative (and that the empty substitution indeed acts as a unit element). Further, clearly the multiplication respects the grading making $(dWHA, m, e)$ a connected graded algebra.

*Comultiplication.* To write down the comultiplication of $dWHA$ a preliminary definition is needed. Let

$$\alpha = [a_1, \cdots, a_m]$$

be a word over an alphabeth $\mathcal{X}$. A good cut of $\alpha$ is a cut $[a_1, \cdots, a_r] \otimes [a_{r+1}, \cdots a_m]$ such that $\mathrm{supp}([a_1, \cdots, a_r]) \cap \mathrm{supp}([a_{r+1}, \cdots, a_m]) = \emptyset$. The two trivial cuts are always good. For example the good cuts of $[x_2, x_3, x_2, x_4, x_1]$ are $1 \otimes [x_2, x_3, x_2, x_4, x_1]$, $[x_2, x_3, x_2] \otimes [x_4, x_1]$, $[x_2, x_3, x_2, x_4] \otimes [x_1]$ and $[x_2, x_3, x_2, x_4, x_1] \otimes 1$ (where as usual 1 is short for [ ]). A subword of $\alpha$ is a word $[a_{i_1}, a_{i_2}, \cdots, a_{i_r}]$ with $i_1 < i_2 < \cdots < i_r$.

The comultiplication of $dWHA$ is now

$$\mu_{dWHA}(p) = \sum_{\sigma_1 \ \sigma_2} \begin{pmatrix} \rho^{-1}(\sigma_1) \\ \sigma_1 \end{pmatrix} \otimes \begin{pmatrix} \rho^{-1}(\sigma_2) \\ \sigma_2 \end{pmatrix} \tag{6.6}$$

where the sum is over all good cuts $\sigma_1 \quad \sigma_2 = \sigma$ of the word $\sigma$. This respects the grading. For instance the comultiplication of the substitution

$$p = \begin{pmatrix} [x_1, x_2, x_1, x_3, x_3, x_1, x_4, x_1, x_4] \\ [x_2, x_3, x_2, x_4, x_1] \end{pmatrix}$$

is

$$\mu_{dWHA}(p) = 1 \quad p + \begin{pmatrix} [x_2, x_3, x_3] \\ [x_2, x_3, x_2] \end{pmatrix} \quad \begin{pmatrix} [x_1, x_1, x_1, x_4, x_1, x_4] \\ [x_4, x_1] \end{pmatrix}$$
$$+ \begin{pmatrix} [x_2, x_3, x_3, x_4, x_4] \\ [x_2, x_3, x_2, x_4] \end{pmatrix} \quad \begin{pmatrix} [x_1, x_1, x_1, x_1] \\ [x_1] \end{pmatrix} + p \quad 1$$

*Counit*. The counit is given by $\varepsilon(p) = 0$ if $\deg(p) > 0$ and $\varepsilon$ takes the value 1 on the empty substitution.

It is now easy to check that $(dWHA, \mu, \varepsilon)$ is a coassociative connected graded coalgebra.

6.7. **Theorem**. The Hopf property holds; i.e $(dWHA, m, \mu, e, \varepsilon)$ is a connected graded bialgebra and hence (because it is connected graded) there is also an antipode, making it a Hopf algebra.

**Proof**. To prove the Hopf property it is needed to prove the commutativity of the following diagram

$$dWHA^{\ 2} \longrightarrow dWHA^{\ 4}$$
$$\downarrow m \qquad\qquad \downarrow id \quad tw \quad id$$
$$dWHA \qquad\qquad dWHA^{\ 4}$$
$$\downarrow \mu \qquad\qquad \downarrow m \quad m$$
$$dWHA^{\ 2} \quad = \quad dWHA^{\ 2}$$

So let

$$p = \begin{pmatrix} \rho \\ \sigma \end{pmatrix} \quad \text{and} \quad p' = \begin{pmatrix} \rho' \\ \sigma' \end{pmatrix}$$

be two substitutions. Their product is

$$m(p \quad p') = \begin{pmatrix} \rho \quad \rho' \\ \sigma \times_{sh} \sigma' \end{pmatrix}$$

Now consider a cut $\gamma = \gamma_1 \quad \gamma_2$ of a shuffle $\gamma$ of $\sigma$ and $\sigma'$. If this is a good cut it induces good cuts of $\sigma$ and $\sigma'$, say $\sigma = \sigma_1 \quad \sigma_2$ and $\sigma' = \sigma'_1 \quad \sigma'_2$ because $\mathrm{supp}(\sigma) \cap \mathrm{supp}(\sigma') = \quad$. Indeed $\sigma_1$ is the prefix of $\sigma$ consisting of all letters of $\sigma$ that occur in the prefix $\gamma_1$ of $\gamma$. (Note that these letters are recognizable because of the support condition and that they form a prefix (not just a subword) because in a shuffle the letters of each of the two factors occur in their original order.) Also $\sigma_2$ is the suffix of $\sigma$ consisting of the letters of $\sigma$ that occur in $\gamma_2$. Similarly one finds a good cut $\sigma' = \sigma'_1 \quad \sigma'_2$ of $\sigma'$. Moreover $\gamma_1$ is a shuffle of $\sigma_1$ and $\sigma'_1$ and $\gamma_2$ is a shuffle of $\sigma_2$ and $\sigma'_2$.

Inversely let $\sigma = \sigma_1 \otimes \sigma_2$ and $\sigma' = \sigma'_1 \otimes \sigma'_2$ be two good cuts, let $\gamma_1$ be any shuffle of $\sigma_1, \sigma'_1$ and $\gamma_2$ a shuffle of $\sigma_2, \sigma'_2$, then $\gamma_1 \otimes \gamma_2$ is a shuffle of $\sigma, \sigma'$ and all shuffles of $\sigma, \sigma'$ are obtained this way. Thus the result of applying $\mu \circ m$ to $p \otimes p'$ is the sum

$$\sum \begin{pmatrix} p^{-1}(\sigma_1) \otimes p'^{-1}(\sigma'_1) \\ \sigma_1 \times_{sh} \sigma'_1 \end{pmatrix} \otimes \begin{pmatrix} p^{-1}(\sigma_2) \otimes p'^{-1}(\sigma'_2) \\ \sigma_2 \times_{sh} \sigma'_2 \end{pmatrix} \tag{6.8}$$

where the sum is over all good cuts $\sigma = \sigma_1 \otimes \sigma_2$, $\sigma' = \sigma'_1 \otimes \sigma'_2$. Here it is also necessary to note that if

$$q = \begin{pmatrix} \rho \otimes \rho' \\ \gamma \end{pmatrix}, \text{ where } \gamma \text{ is a term from the shuffle product of } \sigma, \sigma'$$

and

$$\gamma = \gamma_1 \otimes \gamma_2 \text{ is a good cut with } \gamma_1 \text{ a shuffle of } \sigma_1, \sigma'_1;\ \gamma_2 \text{ a shuffle of } \sigma_2, \sigma'_2$$

then

$$q^{-1}(\gamma_1) = p^{-1}(\sigma_1) \otimes p'^{-1}(\sigma'_1),\ \ q^{-1}(\gamma_2) = p^{-1}(\sigma_2) \otimes p'^{-1}(\sigma'_2)$$

(This uses the support condition again.) However the sum (6.8) is precisely what one gets by applying $(m \otimes m) \circ (id \otimes tw \otimes id) \circ (\mu \otimes \mu)$ to $p \otimes p'$. This proves the theorem (modulo the trivial verifications concerning unit and counit). □

There is a natural nondegenerate non positive definite inner product on $dWHA$ defined by

$$\left\langle \begin{pmatrix} \rho \\ \sigma \end{pmatrix}, \begin{pmatrix} \rho' \\ \sigma' \end{pmatrix} \right\rangle = \begin{cases} 1 & \text{if } \rho = \sigma',\ \rho' = \sigma \\ 0 & \text{otherwise} \end{cases} \tag{6.9}$$

where of course the two substitutions must be so written that $\operatorname{supp}(\rho) = \operatorname{supp}(\sigma) = \operatorname{supp}(\rho') = \operatorname{supp}(\sigma')$. So the more careful statement is that

$$\left\langle \begin{pmatrix} \rho \\ \sigma \end{pmatrix}, \begin{pmatrix} \rho' \\ \sigma' \end{pmatrix} \right\rangle = 1$$

iff there is a substitution of variables (of the variables in $\rho'$ and $\sigma'$) such that $\rho = \sigma'$, $\rho' = \sigma$, and otherwise this inner product is zero.

6.10. **Theorem**. The Hopf algebra $dWHA$ is selfdual with respect to the inner product (6.9).

**Proof**. What needs to be shown is that

$$\left\langle \begin{pmatrix} \rho \\ \sigma \end{pmatrix} \otimes \begin{pmatrix} \rho' \\ \sigma' \end{pmatrix}, \sum \begin{pmatrix} p''^{-1}(\sigma''_1) \\ \sigma''_1 \end{pmatrix} \otimes \begin{pmatrix} p''^{-1}(\sigma''_2) \\ \sigma''_2 \end{pmatrix} \right\rangle = \left\langle \begin{pmatrix} \rho \otimes \rho' \\ \sigma \times_{sh} \sigma' \end{pmatrix}, \begin{pmatrix} \rho'' \\ \sigma'' \end{pmatrix} \right\rangle \tag{6.11}$$

where on the left hand side the sum is over all good cuts $\sigma'' = \sigma''_1 \otimes \sigma''_2$ and where it must be the case (by the definition of the multiplication) that $\mathrm{supp}(\rho)\cap\mathrm{supp}(\rho')=\emptyset$. Now on the right hand there can be a summand that is nonzero only if $\sigma''=\rho*\rho'$ and then because $\rho$ and $\rho'$ have disjoint supports this is a good cut. Moreover $\rho''$ must be a shuffle of $\sigma,\sigma'$. If both these conditions hold the right hand side is 1, otherwise it is zero. Now on the left hand side there can be at most one good cut which yields a term that is nonzero, viz the one for which $\rho=\sigma_1''$ and $\rho'=\sigma_2''$. Further $\rho''$ is made up of the two complimentary subwords $p''^{-1}(\sigma_1'')$ and $p''^{-1}(\sigma_2'')$ with disjoint supports, that is it is a shuffle of these two words. Now for the left hand side to have a term equal to one (there can be only one at most) it must first of all be the case that $\rho=\sigma_1''$ and $\rho'=\sigma_2''$ for some good cut $\sigma''=\sigma_1''*\sigma_2''$ so that also $\sigma''=\sigma_1''*\sigma_2''$ (and there can be at most one good cut like that). It must also be the case that $\sigma=p''^{-1}(\sigma_1'')$ and $\sigma'=p''^{-1}(\sigma_2'')$ making $\rho''$ a shuffle of $\sigma,\sigma'$. Thus if the left hand side of (6.11) is nonzero it is equal to one and then the right hand side is also equal to 1.

Inversely let the righthand side be equal to 1. Then $\sigma''=\rho*\rho'$ and $\rho''$ must be a shuffle of $\sigma,\sigma'$. Let $\sigma_1''$ be the maximal subword of $\sigma''$ with the same support as $\sigma$ as a subword of $\rho''$. Then $\mathrm{supp}(\sigma_1'')=\mathrm{supp}(\sigma)=\mathrm{supp}(\rho)$ and as $\sigma''=\rho*\rho'$ and $\mathrm{supp}(\rho)=\mathrm{supp}(\rho')=\emptyset$ it follows that $\sigma_1''=\rho$ is actually a prefix of $\sigma''$ and that the left hand side of (6.11) is also equal to one. □

## 7. *MPR* as a sub Hopf algebra of *dWHA*.

Let $\tau$ be a permutation word of length $n$. The permutation words over $\mathbf{N}$ are precisely the words of height equal to length and no multiplicities. To $\tau$ associate a substitution as follows

$$\varphi:\ \tau=[t_1,t_2,\cdots,t_n]\ \mapsto p(\tau)=\begin{pmatrix}[x_1,x_2,\cdots,x_n]\\ [x_{t_1},x_{t_2},\cdots,x_{t_n}]\end{pmatrix}\in dWHA \tag{7.1}$$

It is easy to characterize the substitutions that arise this way. They are precisely the substitutions for which both the top word and the bottom word have no multiplcities.

More generally if $\tau=[t_1,t_2,\cdots,t_n]$ is any word over the natural numbers with no multiplicities there is also a permutation substitution attached to it. Indeed, let the support of such a $\tau$ be $\{a_1,a_2,\cdots,a_n\}$, $a_1<a_2<\cdots<a_n$ (same $n$ because of the no multiplicities condition). Then the permutation substitution associated to $\tau$ is

$$\begin{pmatrix}[x_{a_1},x_{a_2},\cdots,x_{a_n}]\\ [x_{t_1},x_{t_2},\cdots,x_{t_n}]\end{pmatrix}$$

which, as a substitution, is the same as $p(\mathrm{st}(\tau))$ as defined by (7.1), where st is the standardization map of section 3, i.e. $\mathrm{st}(\tau)=[\psi(t_1),\psi(t_2),\cdots,\psi(t_n)]$ where $\psi$ is the unique strictly monotone map $\mathrm{supp}(\tau)\longrightarrow\{1,2,\cdots,n\}$. This is how standardization of permutation words appears in this business. Of course also in the world of permutations it is quite customary to identify permutations defined by $\tau$ and $\mathrm{st}(\tau)$. (The underlying alphabeth does not really matter.)

7.2. **Theorem**. The imbedding defined by (7.1) is a monomorphism of connected graded Hopf algebras $MPR\longrightarrow dWHA$.

Here *MPR* has the graded Hopf algebra structure defined in section 3, see especially (3.2), (3.3), and *dWHA* has the Hopf algebra structure described in section 6, see especially formulas (6.3)-(6.4).

**Proof**. That things go well for the grading and for the units and counits is immediate. So let

$$\sigma = [s_1,\cdots,s_m] \quad \text{and} \quad \tau = [t_1,\cdots,t_n]$$

and consider their associated substitutions according to (7.1). The product in $dWHA$ of these substitutions is

$$\begin{pmatrix} [x_1,\cdots,x_m] \quad [y_1,\cdots,y_n] \\ [x_{s_1},\cdots,x_{s_m}] \times_{sh} [y_{t_1},\cdots,y_{t_n}] \end{pmatrix} = \begin{pmatrix} [x_1,\cdots,x_m,x_{m+1},\cdots,x_{m+n}] \\ [x_{s_1},\cdots,x_{s_m}] \times_{sh} [x_{m+t_1},\cdots,x_{m+t_n}] \end{pmatrix}$$

which is the sum of the permutation substitutions corresponding to $\sigma \times_{sh} [m+t_1,\cdots,m+t_n]$. This shows that $\varphi$ preserves multiplication. Now lets look at the comultiplication (in $dWHA$) on $p(\tau)$. This gives

$$\mu(p(\tau)) = \sum_{\text{good cuts}} \begin{pmatrix} p(\tau)^{-1}(\sigma_1) \\ \sigma_1 \end{pmatrix} \quad \begin{pmatrix} p(\tau)^{-1}(\sigma_2) \\ \sigma_2 \end{pmatrix}, \quad p(\tau) = \begin{pmatrix} \rho \\ \sigma \end{pmatrix} \tag{7.3}$$

Now note that because there are no multiplicities all cuts are good cuts, and, again because there are no multiplicities, $\sigma_1$ is a permutation word on the alphabeth formed by the letters in $p(\tau)^{-1}(\sigma_1)$ which is a subword of $[x_1,\cdots,x_n]$ so that the $x$'s in $p(\tau)^{-1}(\sigma_1)$ appear in their natural order. Thus the sum of tensor products of permutation substitutions (7.3) corresponds to

$$\sum_{i=0}^{n} \mathrm{st}([t_1,\cdots,t_i]) \quad \mathrm{st}([t_{i+1},\cdots t_n])$$

proving that $\varphi$ also preserves the comultiplication. □

Things also go well with the duality structure and second multiplication (and, dually, second comultiplication) which are all the obvious restrictions of the analogous structures on $dWHA$.

## 8. The word Hopf algebra *WHA*

This Hopf algebra is a sub Hopf algebra of $dWHA$ that contains $MPR$. It has as basis substitutions of the following form

$$p = \begin{pmatrix} \rho \\ \sigma \end{pmatrix}, \quad \rho = [\underbrace{x_1,x_1,\cdots,x_1}_{r_1},\underbrace{x_2,\cdots x_2}_{r_2},\cdots,\underbrace{x_m,\cdots,x_m}_{r_m}] \tag{8.1}$$

That is besides the support condition on $p$, the top word has the property that if two letters are the same then all the letters between are also equal to these. It is immediate to check that the product of such substitutions is again a sum of substitutions of this kind and also that the coproduct of a substitution of this form is a sum of tensor products of substitutions of this form. (This needs the fact that only 'good cuts' are allowed.) It will save typing (and printing ink and paper) to denote such a word $\rho$ as $[x_1^{r_1},\cdots,x_m^{r_m}]$.

A word of the form (8.1) can be uniquely encoded as a single word over the integers as follows. First let $\alpha = [a_1,a_2,\cdots,a_m]$ be a word over the natural numbers. Let

$$\{a'_1,\cdots,a'_n\} = \mathrm{supp}(\alpha), \quad a'_1 < \cdots < a'_n$$

Then the *WHA* substitution associated to $\alpha$ can be written

$$p(\alpha) = \begin{pmatrix} [x_{a'_1}^{r_1}, \cdots, x_{a'_n}^{r_n}] \\ [x_{a_1}, \cdots, x_{a_m}] \end{pmatrix} \tag{8.2}$$

where $r_1 = a'_1, \cdots, r_i = a'_i - a'_{i-1}, \cdots, r_n = a'_n - a'_{n-1}$. Inversely if $p$ is a substitution of the form (8.1) then the word $\alpha$ associated to it is obtained as follows. Let $\sigma = [x_{i_1}, \cdots, x_{i_t}]$, then $\alpha(p) = [a_1, \cdots, a_t]$ with $a_j = r_1 + \cdots r_{i_j}$.

For instance if $\alpha = [3,2,7,2,4]$

$$p(\alpha) = \begin{pmatrix} [x_2^2, x_3, x_4, x_7^3] \\ [x_3, x_2, x_7, x_2, x_4] \end{pmatrix}$$

and inversely if $p$ is just like above, or written more canonically as

$$p = \begin{pmatrix} [x_1^2, x_2, x_3, x_4^3] \\ [x_2, x_1, x_4, x_1, x_3] \end{pmatrix}$$

$r_1 = 2, r_2 = 1, r_3 = 1, r_4 = 3$ and $i_1 = 2, i_2 = 1, i_3 = 4, i_4 = 1, i_5 = 3$ so that

$$\alpha(p) = [r_1 + r_2, r_1, r_1 + r_2 + r_3 + r_4, r_1, r_1 + r_2 + r_3] = [3,2,7,2,4]$$

Note that if the special substitutions of $WHA$ are written as words over the integers their action on *Shuffle* becomes just like the action defined in section 5 near the end.

For example for the $p$ and $\alpha$ under discussion, $\mathrm{ht}(\alpha) = 7$, it acts as zero on all words of length unequal to 7, and on the special words of length 7 of the form $[x_1, x_1, x_2, x_3, x_4, x_4, x_4]$ exactly like in section 5, i.e. by picking out respectively the third, second, seventh, second, and fourth letter. The difference is that under this interpretation $\alpha$ is also zero on words of *Shuffle* of length 7 that are not of the special form $[x_1, x_1, x_2, x_3, x_4, x_4, x_4]$.

It is also obvious that the multiplication of these special substitutions, when written as words over the integers is like the one in section 5. That is

$$m_{WHA}(\alpha \otimes \beta) = [a_1, \cdots, a_m] \times_{sh} [\mathrm{ht}(\alpha) + b_1, \cdots, \mathrm{ht}(\alpha) + b_n]$$

It is trickier to write down a formula for the comultiplication; but a recipe is of course implied by the remarks above. Here is an example with $p$ and $\alpha$ as above. By definition

$$\mu_{WHA}(p) = 1 \otimes p + \begin{pmatrix} [x_2] \\ [x_2] \end{pmatrix} \otimes \begin{pmatrix} [x_1^2, x_3, x_4^3] \\ [x_1, x_4, x_1, x_3] \end{pmatrix} + \begin{pmatrix} [x_1^2, x_2, x_4^3] \\ [x_2, x_1, x_4, x_1] \end{pmatrix} \otimes \begin{pmatrix} [x_3] \\ [x_3] \end{pmatrix} + p \otimes 1$$

which translated to the level of words over the integers yields, with $\alpha = [3,2,7,2,4]$

$$\mu_{WHA}(\alpha) = 1 \otimes \alpha + [1] \otimes [2,6,2,3] + [3,2,6,2] \otimes [1] + \alpha \otimes 1$$

It is perhaps worth noting that with the interpretation given in this section of words acting as endomorphisms of *Shuffle* the dodgy procedure of section 5 for defining a comultiplication actually works; as it does for $dWHA$.

Michiel Hazewinkel
CWI, Amsterdam
<mich@cwi.nl>